\documentclass[12pt,a4paper]{article}
\usepackage{}
\usepackage{indentfirst}
\usepackage{amsfonts}
\usepackage{amssymb}
\usepackage{mathrsfs}
\usepackage{amsmath}
\usepackage{amsthm}
\usepackage{enumerate}
\usepackage{cite}
\usepackage[all]{xy}
\allowdisplaybreaks
\newtheorem{defn}{Definition}[section]
\newtheorem{thm}[defn]{Theorem}
\newtheorem{lem}[defn]{Lemma}
\newtheorem{prop}[defn]{Proposition}
\newtheorem{cor}[defn]{Corollary}
\newtheorem{eg}[defn]{Example}
\newtheorem{re}[defn]{Remark}
\newcommand{\bdefn}{\begin{defn}}
\newcommand{\edefn}{\end{defn}}
\newcommand{\bthm}{\begin{thm}}
\newcommand{\ethm}{\end{thm}}
\newcommand{\blem}{\begin{lem}}
\newcommand{\elem}{\end{lem}}
\newcommand{\bprop}{\begin{prop}}
\newcommand{\eprop}{\end{prop}}
\newcommand{\bcor}{\begin{cor}}
\newcommand{\ecor}{\end{cor}}
\newcommand{\beg}{\begin{eg}}
\newcommand{\eeg}{\end{eg}}
\newcommand{\bre}{\begin{re}}
\newcommand{\ere}{\end{re}}
\newcommand{\bpf}{\begin{proof}}
\newcommand{\epf}{\end{proof}}

\newcommand{\benu}{\begin{enumerate}}
\newcommand{\eenu}{\end{enumerate}}
\newcommand{\bc}{\begin{center}}
\newcommand{\ec}{\end{center}}
\newcommand{\bea}{\begin{eqnarray}}
\newcommand{\eea}{\end{eqnarray}}
\newcommand{\Bea}{\begin{eqnarray*}}
\newcommand{\Eea}{\end{eqnarray*}}
\newcommand{\beq}{\begin{equation}}
\newcommand{\eeq}{\end{equation}}
\newcommand{\Beq}{\begin{equation*}}
\newcommand{\Eeq}{\end{equation*}}
\newcommand{\bspl}{\begin{split}}
\newcommand{\espl}{\end{split}}
\newcommand\relphantom[1]{\mathrel{\phantom{#1}}}

\numberwithin{equation}{section}

\begin{document}

\title{\bf  Cohomology and deformations of 3-dimensional Heisenberg Hom-Lie superalgebras}
\author{\normalsize \bf Junxia Zhu, Liangyun Chen}
\date{\small{  School of Mathematics and Statistics, Northeast Normal University,\\ Changchun,  130024,  CHINA}}
\maketitle

\begin{abstract}
In this paper, we study Hom-Lie superalgebras of Heisenberg type. For 3-dimensional Heisenberg Hom-Lie superalgebras, we describe their Hom-Lie super structures, compute the cohomology spaces and characterize their infinitesimal deformations.

\end{abstract}

\footnote[0]{Corresponding author(L. Chen): chenly640@nenu.edu.cn.}
\footnote[0]{Supported by  NNSF of China (No. 117701069), NSF of Jilin province(No. 20170101048JC) and  the project of jilin province department of education (No. JJKH20180005K).  }

\section{Introduction}

In recent years, Hom-Lie algebras and other Hom-algebras are widely studied, motivated initially by instances appeared in Physics literature when looking for quantum deformations of some algebras of vector fields. Hom-Lie superalgebras, as a generalization of Hom-Lie algebras, are introduced in \cite{FA}, \cite{FA1}. Furthermore, the cohomology and deformation theories of Hom-algebras are studied in \cite{MA} \cite{FA2}, \cite{FA3}, \cite{YS} and so on, while the two theories of Hom-Lie superalgebras can be seen in \cite{FA}, \cite{LY}.

We will follow \cite{RP}, \cite{RV} to define Heisenberg Hom-Lie superalgebras, which are a special case of 2-step nilpotent Hom-Lie superalgebras. The main idea of this paper is to characterize the infinitesimal deformations of Heisenberg Hom-Lie superalgebras using cohomology.

The paper proceeds as follows. In Section 2, we recall the definitions of Hom-Lie superalgebras. Section $3$ is dedicated to introduce Heisenberg Hom-Lie superalgebras and classify three-dimensional Heisenberg Hom-Lie superalgebras. In Section $4$, we review the cohomology theory and give the 2-nd cohomology spaces of Heisenberg Hom-Lie superalgebras of dimension three. In the last section, we characterize all the infinitesimal deformations of three-dimensional Heisenberg Hom-Lie superalgebras using cohomology.

\section{Preliminaries}

Let $V$ be a vector superspace over a field $\mathbb{F}$; that is, a $\mathbb{Z}_{2}$-graded vector space with a direct sum decomposition $V=V_{\overline{0}} \oplus V_{\overline{1}}$. The elements of $V_{\overline{j}}$, $j=0,1$, are called homogeneous of parity $j$. The parity of homogenous element $x$ is denoted by $|x|$. Moreover, the superspace $End(V)$ have a natural direct sum decomposition $End(V)=End(V)_{\overline{0}} \oplus End(V)_{\overline{1}}$, where $End(V)_{\overline{j}}=\{ f| f(V_{\overline{i}}) \subseteq V_{\overline{i+j}} \}$, $j =0,1$. Elements of $ End(V)_{\overline{j}}$ are homogeneous of parity $j$.

We review the definition of Hom-Lie superalgebra in \cite{FA}.

\bdefn  \label{def1}
A Hom-Lie superalgebra $\mathfrak{g}=(V,[ \cdot , \cdot ], \alpha )$ is a triple consisting of a superspace $V$ over a field $\mathbb{F}$, an even bilinear map $[ \cdot , \cdot ]:V \times V \longrightarrow V$ and an even superspace homomorphism $ \alpha : V \longrightarrow V$ satisfying
\begin{equation}
[x,y]=-(-1)^{|x||y|}[y,x](skew-supersymmetry), \label{eq1}
\end{equation}
\begin{equation}
\circlearrowleft _{x,y,z}(-1)^{|x||z|}[ \alpha (x),[y,z]]=0(hom-Jacobi \ identity) \label{eq2}
\end{equation}
for all homogenous elements $x,y,z \in V$, where $\circlearrowleft _{x,y,z}$ denotes the cyclic summation over $x,y,z$.
\edefn

We denote $\mathfrak{g}_{\overline{0}}=\mathfrak{g}|_{V_{\overline{0}}}$, $\mathfrak{g}_{\overline{1}}=\mathfrak{g}|_{V_{\overline{1}}}$ and then $\mathfrak{g}=\mathfrak{g}_{\overline{0}} \oplus \mathfrak{g}_{\overline{1}}$. It follows that $\mathfrak{g}$ is a Hom-Lie algebra when $\mathfrak{g}_{\overline{1}}=0$. The classical Lie superalgebra can be obtained when $\alpha=\rm{id}$.

A hom-Lie superalgebra is called multiplicative, if $\alpha([x,y])=[\alpha(x),\alpha(y)],\forall x,y$. It is obvious that the classical Lie superalgebras are a special case of multiplicative hom-Lie superalgebras.

The center of Hom-Lie superalgebra $\mathfrak{g}=(V,[ \cdot , \cdot ], \alpha )$ is defined by
$$Z(\mathfrak{g})=\{ x \in V : [x,y]=0,\ \forall y\in V \}.$$

Two hom-Lie superalgebras $(V,[ \cdot , \cdot ]_{1}, \alpha )$ and $(V,[ \cdot , \cdot ]_{2}, \beta )$ are said to be isomorphic if there exists an even bijective homomorphism $\phi:(V,[ \cdot , \cdot ]_{1}) \rightarrow (V,[ \cdot , \cdot ]_{2})$ satisfying
$$\phi([x,y]_{1})=[\phi(x),\phi(y)]_{2},\ \forall x,y \in V$$
$$\ \phi \circ \alpha= \beta \circ \phi.$$
In particular,  $(V,[ \cdot , \cdot ], \alpha )$ and $(V,[ \cdot , \cdot ], \beta )$ are isomorphic if and only if there exists an even automorphism $\phi$ such that $\beta=\phi \alpha \phi^{-1}$.

Let $V$ be a vector superspace as before. A bilinear form $\mathcal{B}$ on $V$ is called homogeneous of parity $j$ if it satisfies $\mathcal{B}(x,y)=0,\forall x,y \in V, |x| \neq |y| + j$;  skew-supersymmetric if $\mathcal{B}(x,y)=-(-1)^{|x||y|}\mathcal{B}(y,x)$ for all homogenous elements $x,y\in V$; non-degenerate if from $\mathcal{B}(x,y)=0$ for all $x \in V$, it follows that $y=0$.

In this paper, we only discuss multiplicative Hom-Lie superalgebras over the complex field $\mathbb{C}$ and the elements mentioned are homogenous.

\section{Heisenberg Hom-Lie superalgebras}

Let $\mathfrak{g}$ be a finite-dimensional Hom-Lie superalgebra with a 1-dimensional homogenous derived ideal such that $[\mathfrak{g},\mathfrak{g}] \subset Z(\mathfrak{g})$. Let $h\in Z(\mathfrak{g})$ be the homogenous generator of $[\mathfrak{g},\mathfrak{g}] $. Then a homogenous skew-supersymmetric bilinear form $\bar{{\mathcal{B}}}$ can be defined on $\mathfrak{g}$ via $[x,y]=\bar{{\mathcal{B}}}(x,y)h$, $\forall x,y\in \mathfrak{g}$. This induces a homogenous skew-supersymmetric bilinear form $\mathcal{B}$ on $\mathfrak{g}/Z(\mathfrak{g})$ via ${\mathcal{B}}(x+Z(\mathfrak{g}),y+Z(\mathfrak{g}))=\bar{{\mathcal{B}}}(x,y)$.

\bdefn
A Hom-Lie superalgebra $\mathfrak{g}$ is called a Heisenberg Hom-Lie superalgebra if the derived ideal $[\mathfrak{g},\mathfrak{g}]$ is generated by a homogenous element $h \in Z(\mathfrak{g})$ and $\mathcal{B}$ is non-degenerate.
\edefn

From now on, we will also denote a Hom-Lie superalgebra by $(\mathfrak{h}, \alpha)$, where $\mathfrak{h}=(V,[ \cdot , \cdot ]_{\mathfrak{h}})$ is a superalgebra and $\alpha$ is an even linear map. All brackets unmentioned in the following are zero.

Let $\mathfrak{g}=( V,[ \cdot , \cdot ], \alpha )$ be a 3-dimensional Heisenberg Hom-Lie superalgebra with a direct sum decomposition $\mathfrak{g}= \mathfrak{g}_{\overline{0}} \oplus \mathfrak{g}_{\overline{1}}$. Let $h \in Z(\mathfrak{g})$ be the homogenous generator of the derived ideal $[\mathfrak{g},\mathfrak{g}]$, We analyze the cases $h \in \mathfrak{g}_{\overline{0}}$ and $h \in \mathfrak{g}_{\overline{1}}$ separately.

\noindent Case $1$. If $h \in \mathfrak{g}_{\overline{0}}$, we have two subcases:

(1.1) There are $u_{1},u_{2} \in \mathfrak{g}_{\overline{0}}$ such that $\{ u_{1},u_{2},h \}$ is a basis of $\mathfrak{g}$ and $[u_{1},u_{2}]=h$, which implies that $\mathfrak{g}$ is a Hom-Lie algebra.

(1.2) There are $v_{1},v_{2} \in \mathfrak{g}_{\overline{1}}$ such that $\{ h | v_{1},v_{2}\}$ is a basis of $\mathfrak{g}$ and $[v_{1},v_{2}]=h$. Then the Hom-Lie superalgebra will be denoted by $(\mathfrak{h}_{1},\alpha)$.

\noindent Case $2$. If $h \in \mathfrak{g}_{\overline{1}}$, there exist $u \in \mathfrak{g}_{\overline{0}}$, $v \in \mathfrak{g}_{\overline{1}}$ such that $\{ u | v,h\}$ is a basis of $\mathfrak{g}$ and $[u,v]=h$. In this case, we denote the Hom-Lie superalgebra by $(\mathfrak{h}_{2},\alpha)$.

\bthm\label{thm1}
Let $\mathfrak{g}$ be a multiplicative Heisenberg Hom-Lie(non-Lie) superalgebra of dimension three. Then $\mathfrak{g}$ must be isomorphic to one of the following:

$(1)$ $(\mathfrak{h}_{1},\left(
          \begin{array}{ccc}
            \mu_{11}\mu_{22} & 0 & 0 \\
            0 & \mu_{11} & 0 \\
            0 & 0 & \mu_{22} \\
          \end{array}
        \right))
,$ $\mu_{11}\mu_{22} \neq 0$;

$(2)$ $(\mathfrak{h}_{1},\left(
          \begin{array}{ccc}
            \mu_{12}\mu_{21} & 0 & 0 \\
            0 & 0 & \mu_{12} \\
            0 & \mu_{21} & 0 \\
          \end{array}
        \right)),
$ $\mu_{12}\mu_{21} \neq 0$;

$(3)$ $(\mathfrak{h}_{1},\left(
          \begin{array}{ccc}
            0 & 0 & 0 \\
            0 & \mu_{11} & \mu_{12} \\
            0 & 0 & 0 \\
          \end{array}
        \right))
;$

$(4)$ $(\mathfrak{h}_{2},\left(
         \begin{array}{ccc}
           \mu_{0} & 0 & 0 \\
           0 & \mu_{11} & 0 \\
           0 & 0 & \mu_{0} \mu_{11} \\
         \end{array}
       \right));$

$(5)$ $(\mathfrak{h}_{2},\left(
         \begin{array}{ccc}
           \mu_{0} & 0 & 0 \\
           0 & \mu_{11} & 0 \\
           0 & 1 & \mu_{11} \\
         \end{array}
       \right)),
(\mu_{0}-1) \mu_{11} =0,$

\noindent where $\mu_{0},\mu_{ij} \in \mathbb{C}$, $i,j=1,2$.

\ethm

\bpf
We analyze the cases $h \in \mathfrak{g}_{\overline{1}}$ and $h \in \mathfrak{g}_{\overline{0}}$ separately.

\noindent Case $1$. If $h \in \mathfrak{g}_{\overline{1}}$, there exists a basis $\{ u | v,h\}$ of $\mathfrak{g}$ such that $[u,v]=h$. Suppose that $ \alpha = \left(
                                                                                                          \begin{array}{ccc}
                                                                                                            \mu_{0} & 0 & 0 \\
                                                                                                            0 & \mu_{11} & \mu_{12} \\
                                                                                                            0 & \mu_{21} & \mu_{22} \\
                                                                                                          \end{array}
                                                                                                        \right)
$, $\mu_{0},\mu_{ij} \in \mathbb{C}$, $i,j=1,2$.

We have that $\mathfrak{g}$ is multiplicative if and only if $\alpha([e_{i},e_{j}])=[\alpha(e_{i}),\alpha(e_{j})]$ for $i,j=1,2,3$, which implies $\mu_{12}=0$ and $\mu_{22}=\mu_{0}\mu_{11}$. Then we obtain that $ \alpha = \left(
                                                                                                          \begin{array}{ccc}
                                                                                                            \mu_{0} & 0 & 0 \\
                                                                                                            0 & \mu_{11} & 0 \\
                                                                                                            0 & \mu_{21} & \mu_{0}\mu_{11} \\
                                                                                                          \end{array}
                                                                                                        \right).$

(a)If $\mu_{21}=0$, we obtain a Heisenberg Hom-Lie superalgebra
                                                                                                        $$(\mathfrak{h}_{2},\left(
         \begin{array}{ccc}
           \mu_{0} & 0 & 0 \\
           0 & \mu_{11} & 0 \\
           0 & 0 & \mu_{0} \mu_{11} \\
         \end{array}
       \right)).$$

(b)If $\mu_{21} \neq 0$, let
$$ \phi = \left(
                                                                                                          \begin{array}{ccc}
                                                                                                            b_{0} & 0 & 0 \\
                                                                                                            0 & b_{11} & 0 \\
                                                                                                            0 & b_{21} & b_{0}b_{11} \\
                                                                                                          \end{array}
                                                                                                        \right),\
 \phi^{-1} = \left(
                                                                                                          \begin{array}{ccc}
                                                                                                            {b_{0}}^{-1} & 0 & 0 \\
                                                                                                            0 & {b_{11}}^{-1} & 0 \\
                                                                                                            0 & -{b_{0}}^{-1}{b_{11}}^{-2}{b_{21}} & {b_{0}}^{-1}{b_{11}}^{-1} \\
                                                                                                          \end{array}
                                                                                                        \right).$$
Then
\begin{equation*}
\begin{split}
\phi \alpha \phi^{-1} &= \left(
                          \begin{array}{ccc}
                            b_{0} & 0 & 0 \\
                            0 & b_{11} & 0 \\
                            0 & b_{21} & b_{0}b_{11} \\
                          \end{array}
                        \right)
\left(
                          \begin{array}{ccc}
                            \mu_{0} & 0 & 0 \\
                            0 & \mu_{11} & 0 \\
                            0 & \mu_{21} & \mu_{0}\mu_{11} \\
                          \end{array}
                        \right)
\left(
                          \begin{array}{ccc}
                                                                                                            {b_{0}}^{-1} & 0 & 0 \\
                                                                                                            0 & {b_{11}}^{-1} & 0 \\
                                                                                                            0 & -{b_{0}}^{-1}{b_{11}}^{-2}{b_{21}} & {b_{0}}^{-1}{b_{11}}^{-1} \\
                                                                                                          \end{array}
                        \right)\\
&=\left(
   \begin{array}{ccc}
     \mu_{0} & 0 & 0 \\
     0 & \mu_{11} & 0 \\
     0 & (1-\mu_{0})\mu_{11}{b_{11}}^{-1} b_{21}+\mu_{21}b_{0} & \mu_{0}a_{11} \\
   \end{array}
 \right).
\end{split}
\end{equation*}

If $\mu_{0} \neq 1$ and $\mu_{11} \neq 0$, then $b_{21}=-(1-\mu_{0})^{-1}{\mu_{11}}^{-1}\mu_{21}b_{0}{b_{11}}$ yields
$$\phi \alpha \phi^{-1}=\left(
                                                                                                          \begin{array}{ccc}
                                                                                                            \mu_{0} & 0 & 0 \\
                                                                                                            0 & \mu_{11} & 0 \\
                                                                                                            0 & 0 & \mu_{0}\mu_{11} \\
                                                                                                          \end{array}
                                                                                                        \right),$$
which induces a Heisenberg Hom-Lie superalgebra like the one in (a).

Otherwise, i.e.$\mu_{0} = 1$ or $\mu_{11}=0$, then $b_{0}= {\mu_{21}}^{-1}$ yields $\phi \alpha \phi^{-1}=\left(
                                                                                                          \begin{array}{ccc}
                                                                                                            \mu_{0} & 0 & 0 \\
                                                                                                            0 & \mu_{11} & 0 \\
                                                                                                            0 & 1 & \mu_{11} \\
                                                                                                          \end{array}
                                                                                                      \right)$. We can obtain a new Heisenberg Hom-Lie superalgebra

                                                                                                       $$(\mathfrak{h}_{2},\left(
         \begin{array}{ccc}
           \mu_{0} & 0 & 0 \\
           0 & \mu_{11} & 0 \\
           0 & 1 & \mu_{11} \\
         \end{array}
       \right)),\ (\mu_{0}-1) \mu_{11} =0.$$

\noindent Case $2$. If $h \in \mathfrak{g}_{\overline{0}}$, there exist $v_{1},v_{2} \in \mathfrak{g}_{\overline{1}}$ such that $\{ h | v_{1},v_{2}\}$ is a basis of $\mathfrak{g}$ and $[v_{1},v_{2}]=h$. In this case, we can get three Heisenberg Hom-Lie superalgebras:
$$(\mathfrak{h}_{1},\left(
          \begin{array}{ccc}
            \mu_{11}\mu_{22} & 0 & 0 \\
            0 & \mu_{11} & 0 \\
            0 & 0 & \mu_{22} \\
          \end{array}
        \right))
(\mu_{11}\mu_{22} \neq 0), \ (\mathfrak{h}_{1},\left(
          \begin{array}{ccc}
            \mu_{12}\mu_{21} & 0 & 0 \\
            0 & 0 & \mu_{12} \\
            0 & \mu_{21} & 0 \\
          \end{array}
        \right))(\mu_{12}\mu_{21} \neq 0)$$
and      $(\mathfrak{h}_{1},\left(
          \begin{array}{ccc}
            0 & 0 & 0 \\
            0 & \mu_{11} & \mu_{12} \\
            0 & 0 & 0 \\
          \end{array}
        \right)).$\epf

 \section{The adjoint cohomology of Heisenberg Hom-Lie superalgebras}

Let $\mathfrak{g}=(V,[\cdot,\cdot],\alpha)$ be a Hom-Lie superalgebra. Let $x_{1},\cdots,x_{k}$ be $k$ homogeneous elements of $V$ and $(x_{1},\cdots,x_{k}) \in \wedge^{k} V$. Then we denote by $|(x_{1},\cdots,x_{k})|=|x_{1}|+\cdots+|x_{k}|$ the parity of $(x_{1},\cdots,x_{k})$. The set $C^{k}_{\alpha}(\mathfrak{g},\mathfrak{g})$ of $k$-hom-cochains of $\mathfrak{g}=(V,[\cdot,\cdot],\alpha)$ is the set of $k$-linear maps $\varphi:\wedge^{k} V \rightarrow V$ satisfying
\begin{equation} \label{eq2}
\varphi(x_{1},\cdots,x_{i+1},x_{i}\cdots,x_{k})=-(-1)^{|x_{i}||x_{i+1}|}\varphi(x_{1},\cdots,x_{i},x_{i+1}\cdots,x_{k}),
\end{equation}
\begin{equation} \label{eq9}
\alpha(\varphi(x_{1},\cdots, x_{k}))=\varphi(\alpha(x_{1}),\cdots, \alpha(x_{k}))
\end{equation}
for $x_{1},\cdots,x_{k} \in V$, $1 \leqslant i \leqslant k-1$. In particular, $C^{0}_{\alpha}(\mathfrak{g},\mathfrak{g})=\{x\in \mathfrak{g} | \alpha(x)=x\}$. Denote by $|\varphi|$ the parity of $\varphi$ and $|\varphi(x_{1},\cdots, x_{k})|=|(x_{1},\cdots, x_{k})|+|\varphi|$. We immediately get a direct sum decomposition $C^{k}_{\alpha}(\mathfrak{g},\mathfrak{g})=C^{k}_{\alpha}(\mathfrak{g},\mathfrak{g})_{\overline{0}} \oplus C^{k}_{\alpha}(\mathfrak{g},\mathfrak{g})_{\overline{1}}$.

A $k$-coboundary operator $\delta^{k}(\varphi):C^{k}_{\alpha}(\mathfrak{g},\mathfrak{g}) \rightarrow C^{k+1}_{\alpha}(\mathfrak{g},\mathfrak{g})$ is defined by
\begin{equation*}\label{eq6}
\begin{split}
\delta^{k}(\varphi)(x_{0},\cdots,x_{k})=& \sum\limits_{0 \leqslant s < t \leqslant k}(-1)^{t+|x_{t}|(|x_{s+1}|+ \cdots +|x_{t-1}|)} \\
& \times \varphi(\alpha(x_{0}),\cdots,\alpha(x_{s-1}),[x_{s},x_{t}],\alpha(x_{s+1}),\cdots,\hat{x_{t}},\cdots,\alpha(x_{k}))\\
& +\sum\limits_{s=1}^{k}(-1)^{s+|x_{s}|(|\varphi|+|x_{0}|+ \cdots +|x_{s-1}|)}[\alpha^{k-1}(x_{s}),\varphi(x_{0},\cdots,\hat{x_{s}},\cdots,x_{k})],
\end{split}
\end{equation*}
where $\hat{x_{i}}$ means that $x_{i}$ is omitted.

The $k$-cocycles space, $k$-coboundaries space and $k$-th cohomology space are defined as:

(1) $Z^{k}(\mathfrak{g},\mathfrak{g})=ker \delta^{k}$, $Z^{k}(\mathfrak{g},\mathfrak{g})_{\overline{j}}=Z^{k}(\mathfrak{g},\mathfrak{g})\cap C^{k}_{\alpha}(\mathfrak{g},\mathfrak{g})_{\overline{j}},\ j=0,1$;

(2) $B^{k}(\mathfrak{g},\mathfrak{g})=Im \delta^{k-1}$, $B^{k}(\mathfrak{g},\mathfrak{g})_{\overline{j}}=B^{k}(\mathfrak{g},\mathfrak{g})\cap C^{k}_{\alpha}(\mathfrak{g},\mathfrak{g})_{\overline{j}},\ j=0,1$;

(3)$H^{k}(\mathfrak{g},\mathfrak{g})=\frac{ Z^{k}(\mathfrak{g},\mathfrak{g}) }{ B^{k}(\mathfrak{g},\mathfrak{g})}=H^{k}(\mathfrak{g},\mathfrak{g})_{\overline{0}} \oplus H^{k}(\mathfrak{g},\mathfrak{g})_{\overline{1}}$, where $H^{k}(\mathfrak{g},\mathfrak{g})_{\overline{j}}=\frac{ Z^{k}(\mathfrak{g},\mathfrak{g})_{\overline{j}} }{ B^{k}(\mathfrak{g},\mathfrak{g})_{\overline{j}}}, \ j=0,1$.

\bthm\label{thm2}
The cohomology spaces of Heisenberg Hom-Lie superalgebras are:

$(1)$ For $\mathfrak{g}=(\mathfrak{h}_{1},\left(
          \begin{array}{ccc}
            \mu_{11}\mu_{22} & 0 & 0 \\
            0 & \mu_{11} & 0 \\
            0 & 0 & \mu_{22} \\
          \end{array}
        \right))(\mu_{11}\mu_{22}\neq0)$,

\noindent $
H^{1}(\mathfrak{g},\mathfrak{g})= \left \langle \left(
                          \begin{array}{ccc}
                            a_{22}+a_{33} & a_{12}\delta_{\mu_{22},1} & a_{13}\delta_{\mu_{11},1} \\
                            0 & a_{22} & 0 \\
                            0 & 0 & a_{33} \\
                          \end{array}
                        \right)
 \right \rangle,
$

\noindent $
H^{2}(\mathfrak{g},\mathfrak{g})= \left \langle \left(
                                        \begin{array}{cccccc}
                                          0 & 0 & 0 & 0 & -\frac{1}{2}\mu_{22}a_{22}\delta_{\mu_{11},1} & -\frac{1}{2}\mu_{11}a_{34}\delta_{\mu_{22},1} \\
                                          0 & \delta_{\mu_{11},1}a_{22} & 0 & 0 & 0 & 0 \\
                                          0 & 0 & 0 & \delta_{\mu_{22},1}a_{34} & 0 & 0 \\
                                        \end{array}
                                      \right)
 \right \rangle ;
$

 $(2)$ For $\mathfrak{g}=(\mathfrak{h}_{1},\left(
          \begin{array}{ccc}
            \mu_{12}\mu_{21} & 0 & 0 \\
            0 & 0 & \mu_{12} \\
            0 & \mu_{21} & 0 \\
          \end{array}
        \right))(\mu_{12}\mu_{21}\neq0)$,

\noindent $
H^{1}(\mathfrak{g},\mathfrak{g})=
\left \langle \left(
                          \begin{array}{ccc}
                            2a_{22} & a_{12}\delta_{\mu_{12}\mu_{21},1} & a_{12}\delta_{\mu_{12}\mu_{21},1} \\
                            0 & a_{22} & 0 \\
                            0 & 0 & a_{22} \\
                          \end{array}
                        \right)
 \right \rangle,
$

\noindent $
H^{2}(\mathfrak{g},\mathfrak{g})= \left \langle     \left(
                                        \begin{array}{cccccc}
                                          0 & 0 & 0 & 0 & 0 & 0 \\
                                          0 & 0 & \mu_{12}a_{23}\delta_{\mu_{12} \mu_{21},1} & -2{\mu_{12}}^{2}a_{23}\delta_{\mu_{12} \mu_{21},1} & 0 & 0 \\
                                          0 & -2\mu_{21}a_{23}\delta_{\mu_{12} \mu_{21},1} & a_{23}\delta_{\mu_{12} \mu_{21},1} & 0 & 0 & 0 \\
                                        \end{array}
                                      \right)
 \right \rangle;$

 $(3)$ For $\mathfrak{g}=( \mathfrak{h}_{1},\left(
          \begin{array}{ccc}
            0 & 0 & 0 \\
            0 & \mu_{11} & \mu_{12} \\
            0 & 0 & 0 \\
          \end{array}
        \right))$,

\noindent $
H^{1}(\mathfrak{g},\mathfrak{g})= \left \langle \left(
                          \begin{array}{ccc}
                            (a_{22}+a_{33})\delta_{\mu_{12},0} & 0 & a_{13} \\
                            0 & a_{22}\delta_{\mu_{12},0} & 0 \\
                            0 & 0 & a_{33}\delta_{\mu_{12},0} \\
                          \end{array}
                        \right)
 \right \rangle,
$

\noindent $
H^{2}(\mathfrak{g},\mathfrak{g})=$

\noindent \ $\left \langle  \left(
                                        \begin{array}{cccccc}
                                          0 & 0 & 0 & \delta_{\mu_{12},0}a_{14} & 0 & a_{16} \\
                                          0 & a_{22}\delta_{\mu_{11}(\mu_{11}-1),0} & a_{23}\delta_{\mu_{11},0}+\mu_{12}a_{22}\delta_{\mu_{11},1} & a_{24}\delta_{\mu_{11},0}+\epsilon a_{22}\delta_{\mu_{11},1} & 0 & \delta_{\mu_{11},0}a_{26} \\
                                          0 & 0 & 0 & 0 & 0 & 0 \\
                                        \end{array}
                                      \right) \right \rangle$,\\
 where $\epsilon={\mu_{11}}^{-1}{\mu_{12}}^{2}$;

 $(4)$ For $\mathfrak{g}=(\mathfrak{h}_{2},\left(
         \begin{array}{ccc}
           \mu_{0} & 0 & 0 \\
           0 & \mu_{11} & 0 \\
           0 & 0 & \mu_{0} \mu_{11} \\
         \end{array}
       \right)
 ),$

\noindent $
H^{1}(\mathfrak{g},\mathfrak{g})= \left \langle \left(
                          \begin{array}{ccc}
                            a_{11} & 0 & 0 \\
                            a_{21}\delta_{\mu_{0},\mu_{11}} & a_{22} & 0 \\
                            a_{31}\delta_{\mu_{0}(\mu_{11}-1),0} & a_{32}\delta_{(\mu_{0}-1)\mu_{11},0} & a_{11}+a_{22} \\
                          \end{array}
                        \right)
 \right \rangle,
$

\noindent $
H^{2}(\mathfrak{g},\mathfrak{g})=$

$\left \langle \left(
                                        \begin{array}{cccccc}
                                          0 & 0 & 0 & 0 & \frac{1}{2}\mu_{0}a_{22}\delta_{\mu_{11},1} & -\mu_{0}a_{34}\delta_{\mu_{0}\mu_{11},1} \\
                                          0 & a_{22}\delta_{\mu_{11},1} & 0 & 0 & 0 & a_{26}\delta_{({\mu_{0}}^{2}-1)\mu_{11},0} \\
                                          0 & a_{32}\delta_{\mu_{11},0} & a_{33}\delta_{\mu_{11},0} & a_{34}\delta_{\mu_{0} \mu_{11}(\mu_{0}\mu_{11}-1),0} & 0 & a_{36}\delta_{\mu_{0}(\mu_{0}-1)\mu_{11},0} \\
                                        \end{array}
                                      \right)
 \right \rangle ;
$

 $(5)$ For $\mathfrak{g}=(\mathfrak{h}_{2},\left(
         \begin{array}{ccc}
           \mu_{0} & 0 & 0 \\
           0 & \mu_{11} & 0 \\
           0 & 1 & \mu_{11} \\
         \end{array}
       \right)
),( \mu_{0}-1) \mu_{11} =0,$

\noindent $
H^{1}(\mathfrak{g},\mathfrak{g})= \left \langle \left(
                          \begin{array}{ccc}
                            0 & 0 & 0 \\
                            0 & a_{33} & 0 \\
                            a_{31}\delta_{\mu_{0},\mu_{11}} & a_{32}\delta_{(\mu_{0}-1)\mu_{11},0} & a_{33} \\
                          \end{array}
                        \right)
 \right \rangle ,
$

\noindent $
H^{2}(\mathfrak{g},\mathfrak{g})=$

$\left \langle \left(
                                        \begin{array}{cccccc}
                                          0 & a_{12}\delta_{\mu_{0},0}\delta_{\mu_{11},0} & a_{13}\delta_{\mu_{0},0}\delta_{\mu_{11},0} & 0 & 0 & a_{16}\delta_{\mu_{0},0} \\
                                          0 & 2a_{33}\delta_{\mu_{11},1}+a_{34}\delta_{\mu_{11},0} & 0 & 0 & \mu_{0} a_{36} & 0 \\
                                          0 & a_{32}\delta_{\mu_{11}(\mu_{11}-1),0} & a_{33}\delta_{\mu_{11}(\mu_{11}-1),0} & a_{34}\delta_{ \mu_{0},2}\delta_{ \mu_{11},0} & 0 & a_{36} \\
                                        \end{array}
                                      \right)
 \right \rangle .
$

\ethm

\bpf

It is easy to obtain $C^{k}_{\alpha}(\mathfrak{g},\mathfrak{g})$ for $k=1,2$ by Eq.(\ref{eq2}) and Eq.(\ref{eq9}).

Taking $\mathfrak{g}=(\mathfrak{h}_{1},\left(
          \begin{array}{ccc}
            \mu_{11}\mu_{22} & 0 & 0 \\
            0 & \mu_{11} & 0 \\
            0 & 0 & \mu_{22} \\
          \end{array}
        \right),\ \mu_{11}\mu_{22} \neq 0$ for example,

\noindent $C^{1}_{\alpha}(\mathfrak{g},\mathfrak{g})=\left(
               \begin{array}{ccc}
                 a_{11} & 0 & 0 \\
                 0 & a_{22} & a_{23}\delta_{\mu_{11},\mu_{22}} \\
                 0 & a_{32}\delta_{\mu_{11},\mu_{22}} & a_{33} \\
               \end{array}
             \right)
$,

\noindent $C^{2}_{\alpha}(\mathfrak{g},\mathfrak{g})=\left(
               \begin{array}{cccccc}
                 0 & a_{12}\delta_{\mu_{11},\mu_{22}} & a_{13} & a_{14}\delta_{\mu_{11},\mu_{22}} & 0 & 0 \\
                 0 & 0 & 0 & 0 & a_{25}\delta_{\mu_{11} \mu_{22},1} & a_{26}\delta_{{\mu_{22}}^{2},1} \\
                 0 & 0 & 0 & 0 & a_{35}\delta_{{\mu_{11}}^{2},1} & a_{36}\delta_{\mu_{11} \mu_{22},1} \\
               \end{array}
             \right)$.

Let $\varphi_{0}=\left(
                  \begin{array}{ccc}
                    a_{11} & 0 & 0 \\
                    0 & a_{22} & a_{23} \\
                    0 & a_{32} & a_{33} \\
                  \end{array}
                \right)
 \in C^{1}_{\alpha}(\mathfrak{g},\mathfrak{g})_{\overline{0}}$. Then
$$B^{2}(\mathfrak{g},\mathfrak{g})_{\overline{0}}=\{ \delta^{1}\varphi_{0}| \varphi_{0} \in C^{1}_{\alpha}(\mathfrak{g},\mathfrak{g})_{\overline{0}} \}=\left\langle \left(
                       \begin{array}{cccccc}
                         0 & 2a_{32} & a_{22}+a_{33}-a_{11} & 2a_{23} & 0 & 0 \\
                         0 & 0 & 0 & 0 & 0 & 0 \\
                         0 & 0 & 0 & 0 & 0 & 0 \\
                       \end{array}
                     \right) \right \rangle
.
$$
Moreover, we have $ \varphi_{0} \in Z^{1}(\mathfrak{g},\mathfrak{g})_{\overline{0}}$ if and only if $\delta^{1}(\varphi_{0})=0$.

In the same way, we suppose $\varphi_{1}=\left(
                  \begin{array}{ccc}
                    0 & a_{12} & a_{13} \\
                    a_{21} & 0 & 0 \\
                    a_{31} & 0 & 0 \\
                  \end{array}
                \right)
 \in C^{1}(\mathfrak{g},\mathfrak{g})_{\overline{1}}$ and immediately get $Z^{1}(\mathfrak{g},\mathfrak{g})_{\overline{1}}$ and $B^{2}(\mathfrak{g},\mathfrak{g})_{\overline{1}}$.

Now suppose $\alpha=\left(
                             \begin{array}{ccc}
                               \mu_{0} & 0 & 0 \\
                               0 & \mu_{11} & \mu_{12} \\
                               0 & \mu_{21} & \mu_{22} \\
                             \end{array}
                           \right)$ and $\psi_{0}=\left(
                      \begin{array}{cccccc}
                        0 & a_{12} & a_{13} & a_{14} & 0 & 0 \\
                        0 & 0 & 0 & 0 & a_{25} & a_{26} \\
                        0 & 0 & 0 & 0 & a_{35} & b_{36} \\
                      \end{array}
                    \right)
 \in C^{2}_{\alpha}(\mathfrak{g},\mathfrak{g})_{\overline{0}}$($\psi_{1}=\left(
                       \begin{array}{cccccc}
                         0 & 0 & 0 & 0 & a_{15} & a_{16} \\
                         0 & a_{22} & a_{23} & a_{24} & 0 & 0 \\
                         0 & a_{32} & a_{33} & a_{34} & 0 & 0 \\
                       \end{array}
                     \right)
 \in C^{2}_{\alpha}(\mathfrak{g},\mathfrak{g})_{\overline{1}}$). We know that $ \psi_{0} \in Z^{2}_{\alpha}(g,g)_{\overline{0}}$\\
 ($\psi_{1} \in Z^{2}_{\alpha}(g,g)_{\overline{1}}$) if and only if $\delta^{1}(\psi_{0})=0$($\delta^{1}(\psi_{1})=0$).\epf

\section{Infinitesimal deformations of Heisenberg Hom-Lie superalgebras}

Let $\mathfrak{g}=(V,[\cdot,\cdot]_{0},\alpha)$ be a Hom-Lie superalgebra and $\varphi:V \times V \longrightarrow V$ be an even bilinear map commuting with $\alpha$. A bilinear map $[\cdot,\cdot]_{t}=[\cdot,\cdot]_{0}+t\varphi(\cdot,\cdot)$
is called an infinitesimal deformation of $\mathfrak{g}$ if $\varphi$ satisfies
\begin{equation}\label{eq7}
[x,y]_{t}=-[y,x]_{t},
\end{equation}
\begin{equation}
\circlearrowleft _{x,y,z} (-1)^{|x||z|}[ \alpha (x),[y,z]_{t}]_{t}=0 \label{eq8}
\end{equation}
for $x,y,z \in V$. The previous E.q.(\ref{eq7}) imply $\varphi$ is skew-supersymmetric. We denote
$$\varphi \circ \psi (x,y,z)= \circlearrowleft _{x,y,z} (-1)^{|x||z|}\varphi( \alpha (x),\psi(y,z)),$$
and then Eq.(\ref{eq8}) can be denoted by $[\cdot,\cdot]_{t} \circ [\cdot,\cdot]_{t}=0$.

\blem\label{lem1}
Let $\mathfrak{g}=(V,[\cdot,\cdot]_{0},\alpha)$ be a Hom-Lie superalgebra and $[\cdot,\cdot]_{t}=[\cdot,\cdot]_{0}+t\varphi(\cdot,\cdot)$ be an infinitesimal deformation of $g=(V,[\cdot,\cdot]_{0},\alpha)$. Then $\varphi \in Z^{2}(\mathfrak{g},\mathfrak{g})_{\overline{0}}$.
\elem

\bpf
By Eq.(\ref{eq8}), we have
\begin{flalign}.\label{eq10}
\begin{split}
0&=[\cdot,\cdot]_{t} \circ [\cdot,\cdot]_{t}\\
&=\circlearrowleft  _{x,y,z} (-1)^{|x||z|}([ \alpha (x),[y,z]_{t}]_{0}+t\varphi (\alpha (x) ,[y,z]_{t}))\\
&=\circlearrowleft  _{x,y,z} (-1)^{|x||z|} \left[ t([\alpha (x),\varphi(y,z)]_{0}+\varphi(\alpha (x) ,[y,z]_{0}))+t^{2}\varphi(\alpha (x),\varphi(y,z)) \right].
\end{split}&
\end{flalign}
Note that
$$\circlearrowleft  _{x,y,z}(-1)^{|x||z|}([\alpha (x),\varphi(y,z)]_{0}+\varphi(\alpha (x) ,[y,z]_{0}))=(-1)^{|x||z|}\delta^{2}\varphi.$$
Hence $\varphi \in Z^{2}(\mathfrak{g},\mathfrak{g})_{\overline{0}}$.\epf

By Eq.(\ref{eq10}), we can see that $[\cdot,\cdot]_{t}$ is an infinitesimal deformation if and only if $\varphi \circ \varphi=0$.

Let $\mathfrak{g}_{t}=(V,[\cdot,\cdot]_{t},\alpha)$ and $\mathfrak{g}_{t}'=(V,[\cdot,\cdot]_{t}',\alpha')$ be two deformations of $\mathfrak{g}$, where $[\cdot,\cdot]_{t}=[\cdot,\cdot]_{0}+t\varphi(\cdot,\cdot)$ and $[\cdot,\cdot]_{t}'=[\cdot,\cdot]_{0}+t\psi(\cdot,\cdot)$. If there exists a linear automorphism $\Phi_{t}=\rm{id}+t\phi(\phi \in C^{1}_{\alpha}(g,g)_{\overline{0}})$, satisfying
$$\Phi_{t}([x,y]_{t})=[\Phi_{t}(x),\Phi_{t}(y)]_{t}^{'}, \ \forall x,y\in V$$
we call the deformations $\mathfrak{g}_{t}$ and $\mathfrak{g}_{t}^{'}$ are equivalent. It is obvious that $\mathfrak{g}_{t}$ and $\mathfrak{g}_{t}^{'}$ are equivalent if and only if $\varphi_{1}-\psi_{1}\in B^{2}(\mathfrak{g},\mathfrak{g})_{\overline{0}}$. Therefore, the set of infinitesimal deformations of $\mathfrak{g}$ can be parameterized by $H^{2}(\mathfrak{g},\mathfrak{g})_{\overline{0}}$.

A deformation $\mathfrak{g}_{t}$ of Hom-Lie superalgebras $\mathfrak{g}$ is called trivial if  it is equivalent to $\mathfrak{g}$.

\bcor \label{cor1}
All the infinitesimal deformations of the following Heisenberg Hom-Lie superalgebras are trivial:

$(1)$ $(\mathfrak{h}_{1}, \left(
                          \begin{array}{ccc}
                            \mu_{11}\mu_{22} & 0 & 0 \\
                            0 & \mu_{11} & 0 \\
                            0 &  & \mu_{22} \\
                          \end{array}
                        \right)), \
\mu_{11}\mu_{22} \neq 0$;

$(2)$ $(\mathfrak{h}_{1},\left(
                          \begin{array}{ccc}
                            \mu_{12}\mu_{21} & 0 & 0 \\
                            0 & 0 & \mu_{12} \\
                            0 & \mu_{21} & 0 \\
                          \end{array}
                        \right))
 $, $\mu_{12}\mu_{21} \neq 0$;

$(3)$ $(\mathfrak{h}_{2},\left(
                          \begin{array}{ccc}
                            \mu_{0} & 0 & 0 \\
                            0 & \mu_{11} & 0 \\
                            0 & 0 & \mu_{0}\mu_{11} \\
                          \end{array}
                        \right)), \
 \mu_{0}({\mu_{0}}^{2}-1) \mu_{11} \neq 0$.
\ecor
\bpf
All the infinitesimal deformations of a Heisenberg Hom-Lie superalgebras are trivial if and only if $H^{2}(\mathfrak{g},\mathfrak{g})_{\overline{0}}=0$.
\epf

In the following, we discuss the non-trivial infinitesimal deformations of Heisenberg Hom-Lie superalgebras. We will distinguish two separate cases: the ones that are also Lie superalgebras and those are not.

We recall the classification of three-dimensional Lie superalgebras(\cite{WZG}):

\bthm \label{thm3}
Let $L=(V,[\cdot,\cdot])$ be a Lie superalgebras with a direct sum decomposition $V=V_{\overline{0}} \oplus V_{\overline{1}}$, where $dimV_{\overline{0}}=1$ and $dimV_{\overline{1}}=2$. There are $e_{1} \in V_{\overline{0}}$ and $e_{2},e_{3} \in V_{\overline{1}}$ such that $\{e_{1} |  e_{2},e_{3} \}$ is a basis of $V$ . Then $L$ must be isomorphic to one of the following:

$$
\begin{array}{l}
  L_{1}:[e_{1},V_{\overline{1}}]=0,[e_{2},e_{2}]=e_{1},[e_{3},e_{3}]=[e_{2},e_{3}]=0; \\
  L_{2}:[e_{1},V_{\overline{1}}]=0,[e_{2},e_{2}]=[e_{3},e_{3}]=e_{1},[e_{2},e_{3}]=0; \\
  L_{3}^{\lambda}:[e_{1},e_{2}]=e_{2},[e_{1},e_{3}]=\lambda e_{3},[V_{\overline{1}},V_{\overline{1}}]=0; \\
  L_{4}:[e_{1},e_{2}]=e_{2},[e_{1},e_{3}]=e_{2}+e_{3},[V_{\overline{1}},V_{\overline{1}}]=0; \\
  L_{5}:[e_{1},e_{2}]=0,[e_{1},e_{3}]=e_{2},[V_{\overline{1}},V_{\overline{1}}]=0.
\end{array}
$$
\ethm

We construct a new Lie superalgebra ${L_{2}}^{'}=(V,[\cdot,\cdot]_{0})$. Let $\{ h|v_{1},v_{2} \}$ be a basis of $V$ satisfying $[v_{1},v_{2}]_{0}=h$. Ther is an even bijective morphism $\phi :(V,[\cdot,\cdot]_{0})\rightarrow L_{2}$
$$\phi(h)=e_{1},\ \phi(v_{1})=e_{2}+ie_{3},\ \phi(v_{2})=\frac{1}{2}e_{2}-\frac{1}{2}ie_{3}  $$
such that $\phi([x,y]_{0})=[\phi(x),\phi(y)], \forall x ,y \in V.$ Then ${L_{2}}^{'}$ is isomorphic to ${L_{2}}$ and we shall replace ${L_{2}}$ with it in Theorem \ref{thm3}.

\bprop \label{prop1}
A non-trivial infinitesimal deformation of Heisenberg Hom-Lie superalgebra $(\mathfrak{h}_{1},\alpha)$, that is also a Lie superalgebra, is isomorphism to
$$({L_{2}}^{'},\left(
                                   \begin{array}{ccc}
                                     0 & 0 & 0 \\
                                     0 & \mu_{11} & \mu_{12} \\
                                     0 & 0 & 0 \\
                                   \end{array}
                                 \right)).
$$
\eprop

\bpf
Denote by $(\mathfrak{h}_{1},\alpha)=(V,[\cdot,\cdot]_{0},\alpha)$. There is a basis $\{ h|v_{1},v_{2} \}$ such that $[v_{1},v_{2}]_{0}=h$ and others are zero.

If $\alpha=\left(
                                   \begin{array}{ccc}
                                     \mu_{11}\mu_{22} & 0 & 0 \\
                                     0 & \mu_{11} & 0 \\
                                     0 & 0 & \mu_{22} \\
                                   \end{array}
                                 \right)
$($\mu_{11}\mu_{22} \neq 0$) or $\left(
                                   \begin{array}{ccc}
                                     \mu_{12}\mu_{21} & 0 & 0 \\
                                     0 & 0 & \mu_{12} \\
                                     0 & \mu_{21} & 0 \\
                                   \end{array}
                                 \right)
$($\mu_{12}\mu_{21} \neq 0$),  all the infinitesimal deformations of $(\mathfrak{g}_{1},\alpha)$ are trivial.

Consider $\alpha=\left(
                                   \begin{array}{ccc}
                                     0 & 0 & 0 \\
                                     0 & \mu_{11} & \mu_{12} \\
                                     0 & 0 & 0 \\
                                   \end{array}
                                 \right)$. Let $\varphi=\left(
               \begin{array}{cccccc}
                 0 & 0 & 0 & a_{14}\delta_{\mu_{12},0} & 0 & 0 \\
                 0 & 0 & 0 & 0 & 0 & a_{26}\delta_{\mu_{11},0} \\
                 0 & 0 & 0 & 0 & 0 & 0 \\
               \end{array}
             \right)
$ be an even 2-cocycle and $[\cdot,\cdot]_{t}=[\cdot,\cdot]_{0}+t\varphi(\cdot,\cdot)$. Then $\varphi \circ \varphi=0$ and $\mathfrak{g}_{t}=(V,[\cdot,\cdot]_{t},\alpha')$ is an infinitesimal deformation of $(\mathfrak{h}_{1},\alpha)$. Moreover, $\mathfrak{g}_{t}$ is a Lie superalgebra if and only if $[\cdot,\cdot]_{t} \circ [\cdot,\cdot]_{t}=0$, i.e. $a_{26}=0$. All deformations are given in Table $1$.
\epf

\begin{tabular}{ccccc}
\multicolumn{5}{c}{Table 1}\\
\hline
$\mu_{11}$ & $\mu_{12}$ & $[\cdot,\cdot]_{t}$ & base change & Hom-Lie superalgebra \\
\hline
$\neq0$ & 0 &  $\begin{array}{c}
                    [v_{1},v_{2}]_{t}=h \\
                    {[v_{2},v_{2}]_{t}=a_{14}h} \\
                  \end{array}$ &  $\left(
                                      \begin{array}{ccc}
                                        1 & 0 & 0 \\
                                        0 & 1 & \frac{1}{2}a_{14} \\
                                        0 & 0 & 1 \\
                                      \end{array}
                                    \right)
                  $   & $({L_{2}}^{'},\left(
                                   \begin{array}{ccc}
                                     0 & 0 & 0 \\
                                     0 & \mu_{11} & {\mu_{12}}' \\
                                     0 & 0 & 0 \\
                                   \end{array}
                                 \right))
$ \\
  &   &  $a_{14} \neq 0$ &  $a_{14} \neq 0$  &  ${\mu_{12}}'=-\frac{1}{2}a_{14}\mu_{11}$  \\
\hline
\end{tabular}

\bprop \label{prop2}
The non-trivial infinitesimal deformations of $(\mathfrak{h}_{2},\alpha)$, that are also Lie superalgebras, are isomorphism to:

\noindent $(1)$ For $\alpha=\left(
                              \begin{array}{ccc}
                                \mu_{0} & 0 & 0 \\
                                0 & \mu_{11} & 0 \\
                                0 & 0 & \mu_{0} \mu_{11} \\
                              \end{array}
                            \right)(\mu_{0}({\mu_{0}}^{2}-1)\mu_{11}=0)
$,

$(a)$ $\left({L_{3}}^{\lambda},\left(
                            \begin{array}{ccc}
                              \mu_{0} & 0 & 0 \\
                              0 & \mu_{11} & 0 \\
                              0 & 0 & \mu_{11} \\
                            \end{array}
                          \right)
\right)(\mu_{0}-1)\mu_{11} =0,$
\quad
$(b)$ $\left({L_{3}}^{0},\left(
                            \begin{array}{ccc}
                              -1 & 0 & 0 \\
                              0 & \mu_{11} & -\mu_{11} \\
                              0 & 0 & -2\mu_{11} \\
                            \end{array}
                          \right)
\right),$

$(c)$ $\left({L_{3}}^{-1},\left(
                            \begin{array}{ccc}
                              0 & 0 & 0 \\
                              0 & \frac{1}{2} \mu_{11} & \frac{1}{2} \xi \mu_{11} \\
                              0 & \frac{1}{2} \xi^{-1}\mu_{11} & \frac{1}{2} \mu_{11} \\
                            \end{array}
                          \right)
\right);$

\noindent $(2)$ For $\alpha=\left(
                              \begin{array}{ccc}
                                \mu_{0} & 0 & 0 \\
                                0 & \mu_{11} & 0 \\
                                0 & 1 & \mu_{11} \\
                              \end{array}
                            \right)
$,

$(a)$ $\left({L_{3}}^{0},\left(
                            \begin{array}{ccc}
                              0 & 0 & 0 \\
                              0 & 0 & 0 \\
                              0 & 0 & \eta \\
                            \end{array}
                          \right)
\right),$
\quad
$(b)$ $\left(L_{4},\left(
                            \begin{array}{ccc}
                              1 & 0 & 0 \\
                              0 & \mu_{11} & 0 \\
                              0 & 0 & \mu_{11} \\
                            \end{array}
                          \right)
\right),$

$(c)$ $\left({L_{3}}^{\mu_{0}},\left(
                            \begin{array}{ccc}
                              \mu_{0} & 0 & 0 \\
                              0 & 0 & {\mu_{0}}^{-1} \\
                              0 & 0 & 0 \\
                            \end{array}
                          \right)
\right),\mu_{0} \neq 0,1.$
\eprop

\bpf
Denote by $(\mathfrak{h}_{2},\alpha)=(V,[\cdot,\cdot]_{0},\alpha)$. There is a basis $\{ u|v,h \}$ such that $[u,v]_{0}=h$ and others are zero.

$(1)$ Consider $\alpha=\left(
                              \begin{array}{ccc}
                                \mu_{0} & 0 & 0 \\
                                0 & \mu_{11} & 0 \\
                                0 & 0 & \mu_{0} \mu_{11} \\
                              \end{array}
                            \right)$, $\mu_{0}({\mu_{0}}^{2}-1)\mu_{11}=0$.

Let $\varphi=\left(
                                      \begin{array}{ccccc}
                                        0 & 0 & 0 & 0 & 0 \\
                                        0 & 0 & 0 & 0 & a_{26}\delta_{({\mu_{0}}^{2}-1)\mu_{11},0} \\
                                        0 & 0 & 0 & 0 & a_{36}\delta_{{\mu_{0}}({\mu_{0}}-1)\mu_{11},0} \\
                                      \end{array}
                                    \right)$ is an even 2-cocycle. Then $\varphi \circ \varphi=0$
and we obtain an infinitesimal deformation $\mathfrak{g}_{t}=((V,[\cdot,\cdot]_{0},\alpha)')$ is an infinitesimal deformation, where $[\cdot,\cdot]_{t}=[\cdot,\cdot]_{0}+t\varphi(\cdot,\cdot)$. It is easy to see that is also a Lie algebra. We analyze the cases $(a)(\mu_{0}-1)\mu_{11}=0$, $(b)\mu_{0}=-1$ and $(c)\mu_{0}=0$ separately, which are given in Table 2.

\begin{table}
\centering
\begin{tabular}{ccc}
\multicolumn{3}{c}{Table 2}\\
\multicolumn{3}{l}{(a) $(\mu_{0}-1)\mu_{11}=0.$}\\
\hline
$[\cdot,\cdot]_{t}$ & base change & Hom-Lie superalgebra \\
\hline
$\begin{array}{ccc}
     [u,v]_{t}=h \\
     {[u,h]_{t}=a_{26}v+a_{36}h} \\
   \end{array}$ & $\left(
   \begin{array}{ccc}
     {k_{1}}^{-1} & 0 & 0 \\
     0 & \tau{a_{26}}^{-1} & -\tau k_{2} \\
     0 & -\tau{a_{36}}^{-1} & \tau k_{1} \\
   \end{array}
 \right)
$ & $\left({L_{3}}^{\lambda},\left(
                            \begin{array}{ccc}
                              \mu_{0} & 0 & 0 \\
                              0 & \mu_{11} & 0 \\
                              0 & 0 & \mu_{11} \\
                            \end{array}
                          \right)
\right)$ \\
$a_{26}a_{36} \neq 0$ & $\begin{array}{ccc}
{\tau=(k_{1}-k_{2})^{-1}} \\
     k_{1} \neq k_{2},k_{1}k_{2}=-{a_{26}}^{-1} \\
     {k_{1}+k_{2}=-a_{36}{a_{26}}^{-1}} \\
   \end{array}$  & $\begin{array}{ccc}
     \lambda=k_{1}{k_{2}}^{-1}, \lambda \neq 0,-1 \\
     {(\mu_{0}-1)\mu_{11} =0} \\
   \end{array}$ \\
\hline
$\begin{array}{ccc}
     [u,v]_{t}=h \\
     {[u,v]_{t}=a_{36}h} \\
   \end{array}$ & $\left(
   \begin{array}{ccc}
     a_{36} & 0 & 0 \\
     0 & {a_{36}}^{-1} & 1 \\
     0 & {a_{36}}^{-1} & 0 \\
   \end{array}
 \right)
$ & $\left({L_{3}}^{0},\left(
                            \begin{array}{ccc}
                              \mu_{0} & 0 & 0 \\
                              0 & \mu_{11} & 0 \\
                              0 & 0 & \mu_{11} \\
                            \end{array}
                          \right)
\right)$ \\
$a_{36} \neq 0$ &  $a_{36} \neq 0$  & $(\mu_{0}-1)\mu_{11}=0$  \\
\hline
$\begin{array}{ccc}
     [u,v]_{t}=h \\
     {[u,v]_{t}=a_{26}v} \\
   \end{array}$ & $\left(
   \begin{array}{ccc}
     {a_{26}}^{\frac{1}{2}} & 0 & 0 \\
     0 & \kappa & \kappa {a_{26}}^{\frac{3}{2}} \\
     0 & - \kappa {a_{26}}^{\frac{1}{2}} & \kappa a_{26} \\
   \end{array}
 \right)
$ & $\left({L_{3}}^{-1},\left(
                            \begin{array}{ccc}
                              \mu_{0} & 0 & 0 \\
                              0 & \mu_{11} & 0 \\
                              0 & 0 & \mu_{11} \\
                            \end{array}
                          \right)
\right)$ \\
$a_{26} \neq 0$ &  $\begin{array}{c}
                      \kappa=(1+a_{26})^{-1} \\
                      {a_{26} \neq 0}
                    \end{array}$  & $(\mu_{0}-1)\mu_{11}=0$  \\
\hline
\multicolumn{3}{c}{}\\
\multicolumn{3}{l}{(b) $\mu_{0}=-1.$}\\
\hline
$[\cdot,\cdot]_{t}$ & base change & Hom-Lie superalgebra \\
\hline
$\begin{array}{ccc}
     [u,v]_{t}=h \\
     {[u,v]_{t}=a_{36}h} \\
   \end{array}$ & $\left(
   \begin{array}{ccc}
     a_{36} & 0 & 0 \\
     0 & {a_{36}}^{-1} & 1 \\
     0 & {a_{36}}^{-1} & 0 \\
   \end{array}
 \right)
$ & $\left({L_{3}}^{0},\left(
                            \begin{array}{ccc}
                              -1 & 0 & 0 \\
                              0 & -\mu_{11} & \mu_{11} \\
                              0 & 0 & 2\mu_{11} \\
                            \end{array}
                          \right)
\right)$ \\
$a_{36} \neq 0$ &  $a_{36} \neq 0$  &    \\
\hline
\multicolumn{3}{c}{}\\
\multicolumn{3}{l}{(c) $\mu_{0}=0.$}\\
\hline
$[\cdot,\cdot]_{t}$ & base change & Hom-Lie superalgebra \\
\hline
$\begin{array}{ccc}
     [u,v]_{t}=h \\
     {[u,v]_{t}=a_{26}v} \\
   \end{array}$ & $\left(
   \begin{array}{ccc}
     {a_{26}}^{\frac{1}{2}} & 0 & 0 \\
     0 & \rho & \rho {a_{26}}^{\frac{3}{2}} \\
     0 & - \rho {a_{26}}^{\frac{1}{2}} & \rho a_{26} \\
   \end{array}
 \right)
$ & $\left({L_{3}}^{0},\left(
                            \begin{array}{ccc}
                              0 & 0 & 0 \\
                              0 & \frac{1}{2}\mu_{11} & \frac{1}{2}\xi \mu_{11} \\
                              0 & \frac{1}{2}\xi^{-1}\mu_{11} & \frac{1}{2}\mu_{11} \\
                            \end{array}
                          \right)
\right)$ \\
 $a_{26} \neq 0$ &  $\begin{array}{c}
                      \rho=(1+a_{26})^{-1} \\
                      {a_{26} \neq 0}
                    \end{array}$
   &  $\xi=-{a_{26}}^{-\frac{1}{2}}$  \\
\hline
\end{tabular}
\end{table}

$(2)$ Consider $\alpha=\left(
                              \begin{array}{ccc}
                                \mu_{0} & 0 & 0 \\
                                0 & \mu_{11} & 0 \\
                                0 & 1 & \mu_{11} \\
                              \end{array}
                            \right), \ (\mu_{0}-1)\mu_{11}=0$.

Let $\varphi=\left(
                                      \begin{array}{cccccc}
                                        0 & a_{12}\delta_{\mu_{0},0}\delta_{\mu_{11},0} & a_{13}\delta_{\mu_{0},0}\delta_{\mu_{11},0} & 0 & 0 & 0 \\
                                        0 & 0 & 0 & 0 & \mu_{0} a_{36} & 0 \\
                                        0 & 0 & 0 & 0 & 0 & a_{36} \\
                                      \end{array}
                                    \right)$ be an even 2-cocycle. Then $\mathfrak{g}_{t}=((V,[\cdot,\cdot]_{0},\alpha)')$ is an infinitesimal deformation, where $[\cdot,\cdot]_{t}=[\cdot,\cdot]_{0}+t\varphi(\cdot,\cdot)$, if and only if $\varphi \circ \varphi=0$. We analyze the cases $(a)\mu_{0}=\mu_{11}=0$, $(b)\mu_{0}=1$ and $(c)\mu_{0}\neq 0,1, \ \mu_{11}=0$ separately.

For case (a), $\varphi \circ \varphi=0$ implies $a_{12}=a_{13}=0$ or $a_{36}=0$. Furthermore, if $a_{12}=a_{13}=0$, the deformation $\mathfrak{g}_{t}$ is also a Lie superalgebra.

For cases (b) and (c), $\mathfrak{g}_{t}$ is an infinitesimal deformation for all $\varphi$. The deformation is also a Lie superalgebra if $a_{12}=a_{13}=0$.

The deformations of $(a),\ (b)$ and $(c)$ are given in Table 3.
\begin{table}
\centering
\begin{tabular}{ccc}
\multicolumn{3}{c}{Table 3}\\
\multicolumn{3}{l}{(a) $\mu_{0}=\mu_{11}=0.$}\\
\hline
$[\cdot,\cdot]_{t}$ & base change & Hom-Lie superalgebra \\
\hline
$\begin{array}{ccc}
     [u,v]_{t}=h \\
     {[u,v]_{t}=a_{36}h} \\
   \end{array}$ & $\left(
   \begin{array}{ccc}
     a_{36} & 0 & 0 \\
     0 & 0 & {a_{36}}^{-1} \\
     0 & 1 & {a_{36}}^{-1} \\
   \end{array}
 \right)
$ & $\left({L_{3}}^{0},\left(
                            \begin{array}{ccc}
                              0 & 0 & 0 \\
                              0 & 0 & 0 \\
                              0 & 0 & \mu_{22} \\
                            \end{array}
                          \right)
\right)$ \\
$a_{36} \neq 0$ &  $a_{36} \neq 0$  &  $\mu_{22}=a_{36}$  \\
\hline
\multicolumn{3}{c}{}\\
\multicolumn{3}{l}{(b) $\mu_{0} = 1.$}\\
\hline
$[\cdot,\cdot]_{t}$ & base change & Hom-Lie superalgebra \\
\hline
$\begin{array}{ccc}
     [u,v]_{t}=a_{36}v+h \\
     {[u,v]_{t}=a_{36}h} \\
   \end{array}$ & $\left(
   \begin{array}{ccc}
     a_{36} & 0 & 0 \\
     0 & 0 & a_{36} \\
     0 & 1 & 0 \\
   \end{array}
 \right)
$ & $\left(L_{4},\left(
                            \begin{array}{ccc}
                              1 & 0 & 0 \\
                              0 & \mu_{11} & 0 \\
                              0 & 0 & \mu_{11} \\
                            \end{array}
                          \right)
\right)$ \\
$a_{36} \neq 0$ &  $a_{36} \neq 0$  &    \\
\hline
\multicolumn{3}{c}{}\\
\multicolumn{3}{l}{(c) $\mu_{0}\neq 0,1,\ \mu_{11}= 0.$}\\
\hline
$[\cdot,\cdot]_{t}$ & base change & Hom-Lie superalgebra \\
$\begin{array}{ccc}
     [u,v]_{t}=\mu_{0} a_{36}v+h \\
     {[u,v]_{t}=a_{36}h} \\
   \end{array}$ & $\left(
   \begin{array}{ccc}
     a_{36} & 0 & 0 \\
     0 & (\mu_{0}-1)a_{36} & {\mu_{0}}^{-1} \\
     0 & 0 & 0 \\
   \end{array}
 \right)
$ & $\left({L_{3}}^{\lambda},\left(
                            \begin{array}{ccc}
                              \mu_{0} & 0 & 0 \\
                              0 & 0 & {\mu_{0}}^{-{1}} \\
                              0 & 0 & 0 \\
                            \end{array}
                          \right)
\right)$ \\
 $a_{36} \neq 0$ &  $a_{36} \neq 0$  &  $\lambda \neq 0,1$  \\
\hline
\end{tabular}
\end{table}\epf

Proposition \ref{prop1} and Proposition \ref{prop2} give the infinitesimal deformations of Heisenberg Hom-Lie superalgebras that are also Lie superalgebras. Before discussing the rest deformations, we will recall some multiplicative Hom-Lie superalgebras and those can be find in the classification of multiplicative Hom-Lie superalgebras of \cite{WCY}. Let $V$ be a  superspace with a direct sum decomposition $V=V_{\overline{0}} \oplus V_{\overline{1}}$, $[\cdot,\cdot]$ be an even bilinear map and $\sigma$ be an even linear map on $V$. Let $\{e_{1} | e_{2},e_{3}\}$ be a basis of $V$.  The following are three Hom-Lie superalgebras on $V$:

\noindent $L_{1,2}^{43,a}:[e_{1},e_{2}]=0,[e_{1},e_{3}]=\beta e_{2}, [e_{2},e_{2}]=0,[e_{2},e_{3}]=0,[e_{3},e_{3}]=\gamma e_{1}, \sigma=\left(
                                                                                 \begin{array}{ccc}
                                                                                   0 & 0 & 0 \\
                                                                                   0 & 0 & a \\
                                                                                   0 & 0 & 0 \\
                                                                                 \end{array}
                                                                               \right)$ $(a \neq 0),$

\noindent $L_{1,2}^{45,a}:[e_{1},e_{2}]=0,[e_{1},e_{3}]=\beta e_{2}, [e_{2},e_{2}]=0,[e_{2},e_{3}]=\nu e_{1},[e_{3},e_{3}]=\gamma e_{1}, \sigma=\left(
                                                                                 \begin{array}{ccc}
                                                                                   0 & 0 & 0 \\
                                                                                   0 & 0 & a \\
                                                                                   0 & 0 & 0 \\
                                                                                 \end{array}
                                                                               \right)$ $(a \neq 0),$

\noindent $L_{1,2}^{46,a,b}:[e_{1},e_{2}]=0,[e_{1},e_{3}]=\beta e_{2}, [e_{2},e_{2}]=0,[e_{2},e_{3}]=\mu e_{1},[e_{3},e_{3}]=0, \sigma=\left(
                                                                                 \begin{array}{ccc}
                                                                                   0 & 0 & 0 \\
                                                                                   0 & 0 & a \\
                                                                                   0 & 0 & b \\
                                                                                 \end{array}
                                                                               \right)$ $ (a^{2}+b^{2} \neq 0)$.

The following Proposition \ref{prop3} and Proposition \ref{prop4} will characterize the infinitesimal deformations of Heisenberg Hom-Lie superalgebras that are not Lie superalgebras.

\bprop \label{prop3}
An infinitesimal deformation of Heisenberg Hom-Lie superalgebra $(\mathfrak{h}_{1},\alpha)$, which is not a Lie superalgebra, is isomorphic to $L_{1,2}^{46,a,0}(a \neq 0)$.
\eprop

\bpf
By the proof of Proposition \ref{prop1}, we know that $(\mathfrak{h}_{1},\alpha)$ have an infinitesimal deformation(not a Lie superalgebra) if and only if
$$\alpha=\left(
                                                                                                                        \begin{array}{ccc}
                                                                                                                          0 & 0 & 0 \\
                                                                                                                          0 & 0 & \mu_{12} \\
                                                                                                                          0 & 0 & 0 \\
                                                                                                                        \end{array}
                                                                                                                      \right)(\mu_{12}\neq 0)
, \ \varphi=\left(
                 \begin{array}{cccccc}
                   0 & 0 & 0 & 0 & 0 & 0 \\
                   0 & 0 & 0 & 0 & 0 & a_{26} \\
                   0 & 0 & 0 & 0 & 0 & 0 \\
                 \end{array}
               \right)
(a_{26} \neq 0).$$
There is a basis $\{h| v_{1}, v_{2} \}$ of $V$ such that $[v_{1},v_{2}]_{0}=h$. Therefore $\mathfrak{g}_{t}=(V,[\cdot,\cdot]_{t},\alpha')$ is an infinitesimal deformation, where $[\cdot,\cdot]_{t}=[\cdot,\cdot]_{0}+t\varphi$ and $[v_{1},v_{2}]_{t}=h$, $[h,v_{2}]_{t}=a_{26}v_{1}(a_{26} \neq 0).$ Then the deformation $\mathfrak{g}_{t}$ is isomorphic to $L_{1,2}^{46,\mu_{12},0}:$

\noindent \ $[e_{1},e_{2}]=0,[e_{1},e_{3}]=a_{26} e_{2}, [e_{2},e_{2}]=0,[e_{2},e_{3}]= e_{1},[e_{3},e_{3}]=0, \sigma=\left(
                                                                                 \begin{array}{ccc}
                                                                                   0 & 0 & 0 \\
                                                                                   0 & 0 & \mu_{12} \\
                                                                                   0 & 0 & 0 \\
                                                                                 \end{array}
                                                                               \right)$.
\epf

\bprop \label{prop4}
An infinitesimal deformation of $(\mathfrak{h}_{2},\alpha)$, which is not a Lie superalgebra, is isomorphic to
$L_{1,2}^{43,1}$, $L_{1,2}^{45,1}$, or $L_{1,2}^{46,1,0}$.
\eprop

\bpf
By the proof of Proposition \ref{prop2}, Heisenberg Hom-Lie superalgebra $(\mathfrak{h}_{2},\alpha)$ have an infinitesimal deformations if and only if $$\alpha=\left(
                                                                                                                        \begin{array}{ccc}
                                                                                                                          0 & 0 & 0 \\
                                                                                                                          0 & 0 & 0 \\
                                                                                                                          0 & 1 & 0 \\
                                                                                                                        \end{array}
                                                                                                                      \right)
, \ \varphi=\left(
                 \begin{array}{cccccc}
                   0 & a_{12} & a_{13} & 0 & 0 & 0 \\
                   0 & 0 & 0 & 0 & 0 & 0 \\
                   0 & 0 & 0 & 0 & 0 & 0 \\
                 \end{array}
               \right).$$
There is a basis $\{ u|v,h \}$ of $V$ such that $[u,v]_{0}=h$. Therefore $\mathfrak{g}_{t}=(V,[\cdot,\cdot]_{t},\alpha')$ is an infinitesimal deformation, where $[\cdot,\cdot]_{t}=[\cdot,\cdot]_{0}+t\varphi$ and $[u,v]_{t}=h$, $[v,v]_{t}=a_{12}u$, $[v,h]_{t}=a_{13}u$. We analyze into three cases.

$(a)$ If $a_{12} \neq 0$ and $a_{13} \neq 0$, the infinitesimal deformation $\mathfrak{g}_{t}$ is isomorphic to

\noindent $L_{1,2}^{45,1}:[e_{1},e_{2}]=0,[e_{1},e_{3}]= e_{2}, [e_{2},e_{2}]=0,[e_{2},e_{3}]=a_{13}e_{1},[e_{3},e_{3}]=a_{12} e_{3}, \sigma=\left(
                                                                                 \begin{array}{ccc}
                                                                                   0 & 0 & 0 \\
                                                                                   0 & 0 & 1 \\
                                                                                   0 & 0 & 0 \\
                                                                                 \end{array}
                                                                               \right)$.

$(b)$ If $a_{12} = 0$ and $a_{13} \neq 0$, $\mathfrak{g}_{t}$ is isomorphic to:

\noindent $L_{1,2}^{46,1,0}:[e_{1},e_{2}]=0,[e_{1},e_{3}]=e_{2}, [e_{2},e_{2}]=0,[e_{2},e_{3}]=a_{13} e_{1},[e_{3},e_{3}]=0, \sigma=\left(
                                                                                 \begin{array}{ccc}
                                                                                   0 & 0 & 0 \\
                                                                                   0 & 0 & 1 \\
                                                                                   0 & 0 & 0 \\
                                                                                 \end{array}
                                                                               \right)$.

$(c)$ If $a_{12} \neq 0$ and $a_{13} = 0$, $\mathfrak{g}_{t}$ is isomorphic to:

\noindent $L_{1,2}^{43,1}:[e_{1},e_{2}]=0,[e_{1},e_{3}]=e_{2}, [e_{2},e_{2}]=0,[e_{2},e_{3}]=0,[e_{3},e_{3}]=a_{12} e_{1}, \sigma=\left(
                                                                                 \begin{array}{ccc}
                                                                                   0 & 0 & 0 \\
                                                                                   0 & 0 & 1 \\
                                                                                   0 & 0 & 0 \\
                                                                                 \end{array}
                                                                               \right)$. \epf

\end{document}